\documentclass[12pt]{article}
\usepackage[cp850]{inputenc}
\usepackage{latexsym,amsfonts}

\newtheorem{Thm}{Theorem}[section]

\newtheorem{Rem}[Thm]{Remark}
\newtheorem{Con}[Thm]{Conjecture}

\author{ Vladimir D. Tonchev\thanks{ Research supported by NSA Grant
H98230-16-1-0011.}\\
 Michigan Technological University\\
 Houghton, Michigan 49931, USA
}
\title{\bf On resolvable Steiner 2-designs and maximal arcs in projective planes}

\begin{document}
\maketitle

\begin{center}
{\it In Honor of Andries Brouwer's 65th Birthday}
\end{center}

\begin{abstract}
A combinatorial characterization of resolvable Steiner 2-$(v,k,1)$ designs
embeddable as maximal arcs in a projective plane of order $(v-k)/(k-1)$
is proved, and a generalization of a conjecture by Andries Brouwer \cite{Br} 
is formulated.
 
\end{abstract}

\vspace{5mm}

\section{Introduction}

We assume familiarity with basic facts and notions from 
combinatorial design theory  \cite{AK},
 \cite{BJL}, \cite{CRC},  \cite{T88}.

Let $D=\{ X, \cal{B} \}$ be a Steiner 2-$(v,k,1)$ design with
point set $X$, collection of blocks $\cal{B}$, and let $v$ be a multiple
of $k$: $v=nk$. Since every point of $X$ is contained in
\[ r=(v-1)/(k-1)=(nk-1)/(k-1) \]
 blocks, it follows that $k-1$ divides $n-1$.
Thus, $n-1=s(k-1)$ for some integer $s\ge 1$, and
\[ v=nk = (sk -s +1)k. \]

A {\it parallel class} (or {\it spread}) is a set of $v/k=n$ 
pairwise disjoint blocks,
and a {\it resolution} of $D$ is a partition of the collection of blocks
$\cal{B}$ into $r=(v-1)/(k-1)=sk+1$ parallel classes.
A design is {\it resolvable} if it admits a resolution.

Any 2-$((sk -s +1)k, k, 1)$ design with $s=1$ is equivalent to
an affine plane of order $k$, and admits exactly one resolution. 
If $s>1$, a resolvable 2-$((sk -s +1)k, k, 1)$ design may admit more than one
resolution.

A property of resolvable Steiner 2-designs having several 
resolutions, that has attracted considerable attention, is 
{\it orthogonality} (see \cite[page 31]{BJL},
 \cite[Section II.7.7]{CRC}, \cite{HV}, 
\cite{LV}, and the references within): 
two resolutions $R_1$, $R_2$, 
\begin{equation}
\label{def}
 R_1 = P^{(1)}_1 \cup P^{(1)}_2 \cup \cdots P^{(1)}_r, \  R_2 =  
 P^{(2)}_1 \cup P^{(2)}_2 \cup \cdots P^{(2)}_r 
\end{equation}
are called {\it orthogonal} if
\[  |P^{(1)}_i \cap P^{(2)}_j | \le 1, \ {\rm for \ all} \ 1\le i, j \le r, \]
that is, every two parallel classes $P^{(1)}_i, \  P^{(2)}_j$,
one from each resolution, share at most one block.

The subject of this paper is a concept which is somewhat similar,
but yet different from orthogonality.
We call two resolutions $R_1$, $R_2$ (\ref{def})
{\em compatible} if they share one parallel class,
\[ P^{(1)}_i = P^{(2)}_j, \]
and
\[ |P^{(1)}_{i'} \cap P^{(1)}_{j'}|\le 1 \]
for  $i' \neq i$ and  $j' \neq j$.

More generally, a set of $m$ resolutions
$ R_{1}, \ldots, R_{m}$ 
 is  compatible 
if every two of these resolutions are  compatible.

Sets of mutually compatible resolutions arise naturally in resolvable
Steiner designs associated with maximal arcs in projective planes
\cite[Section 8.4]{AK}, \cite{Bru}.

In this paper, we prove an upper bound on the maximum number 
of mutually  compatible resolutions, and give a characterization
of the case when this maximum is achieved.

\section{An upper bound on the number of mutually  compatible resolutions}

Suppose that $\cal{P}$ is a projective plane of order $q=sk$. A {\it maximal}
$\{ (sk-s+1)k; k \}$-{\it arc} \cite{Hir}, \cite[p. 558]{JAT}, 
is a set $\cal{A}$ of $(sk-s+1)k$ points
of $\cal{P}$ such that every line  of $\cal{P}$ is ether disjoint from $\cal{A}$
or meets $\cal{A}$ in exactly $k$ points.
The set of lines of $\cal{P}$ which have no points in common with $\cal{A}$
determines a maximal  $\{ (sk-k+1)s; s)$-arc ${\cal{A}}^{\perp}$ in the dual plane
${\cal{P}}^{\perp}$.

A maximal arc with $k=2$ is called a {\it hyperoval} (or oval).
 Maximal arcs, and hyperovals in particular, have been studied in connection
with the construction of projective planes  \cite[Section 8.4]{AK}, \cite{Bru},
\cite{Lam}, \cite{Th}, and partial geometries \cite{JAT}.

Maximal arcs with $1<k < q$ do not exist in any Desarguesian plane of odd order
$q$ \cite{BBM}, 
and are known to exist in every Desarguesian plane of
even order $q$ \cite{Den}, for any $k=2^i, \ k<q$,
  as well as in some non-Desarguesian planes of even order
\cite{HST}, \cite{PRS}.

If $k>1$, the non-empty intersections of
a maximal $\{ (sk-s+1)k; k \}$-arc $\cal{A}$ with the lines of
a projective plane $\cal{P}$ of order $q=sk$
are the blocks of a resolvable 2-$((sk-s+1)k,k,1)$ design $D$.
Similarly, if $s>1$, the corresponding $\{ (sk-k+1)s; s \}$-arc ${\cal{A}}^{\perp}$
in the dual plane is the point set of a resolvable
2-$((sk-k+1)s, s, 1)$ design $D^{\perp}$.
We will refer to $D$ (resp. $D^{\perp}$) as a design embeddable in
$\cal{P}$ (resp. ${\cal{P}}^\perp$) as a maximal arc.
The points of $D^{\perp}$ determine a set of $(sk-k+1)s$ mutually compatible
resolutions of $D$. Respectively, the points of $D$ determine a set
of $(sk-s+1)k$ mutually compatible resolutions of $D^{\perp}$.

\begin{Thm}
\label{th}
Let $S=\{ R_1,\ldots, R_m \}$ be a set of $m$ mutually  compatible resolutions
of a  2-$((sk -s +1)k, k, 1)$ design $D=\{ X, \cal{B} \}$. 
Then
\[ m\le (sk-k+1)s. \]
The equality
\[ m=(sk-k+1)s \]
holds if and only if there exists a projective plane $\cal{P}$ of order $sk$
such that $D$ is embeddable in $\cal{P}$ as a maximal $\{ (sk-s+1)k; k \}$-arc.

\end{Thm}

{\bf Proof}.
It is straightforward to check that if $\cal{P}$
is a projective plane of order $sk$ in which $D$ is embedded as
 a maximal $\{ (sk-s+1)k; k \}$-arc, then the points of 
the corresponding maximal $\{ (sk-k+1)s; s \}$-arc
in the dual plane ${\cal{P}^\perp}$ define a set of $(sk-k+1)s$
mutually  compatible resolutions of $D$.

To prove the converse, we consider the simple incidence structure $\cal{I}$ 
having as "points"
the blocks of $D$, that is, having $\cal{B}$ as a point set, and
having  blocks of size
\[ v/k = n =sk-s+1 \]
being the parallel classes of $D$ which appear in resolutions of $S$.
Let $r_i$ denote the number of blocks of $\cal{I}$ containing the $i$th point of 
 $\cal{I}$.

A point of $D$ is contained in
\[  r=sk+1 \]
blocks, and the total number $b$ of blocks of $D$ is equal to 
\[ b=\frac{vr}{k}=(sk-s+1)(sk+1). \]
Any block of $D$ is disjoint from exactly
\[ b-1 -k(r-1)= s(sk-k+1)(k-1) \]
other blocks of $D$.
It follows that every block of $D$ is contained in at most
\[
 \frac{b-1 -k(r-1)}{\frac{v}{k}-1}=sk-k+1     
\]
parallel classes of $D$ which appear in resolutions from $S$,
that is,
\[ r_i \le sk -k+1. \]
Let $b_{I}$ denote the number of blocks of  $\cal{I}$.
Counting in two ways the incident pairs of a block and a point of 
 $\cal{I}$ gives
\[ b_{I}(\frac{v}{k})=\sum_{i=1}^b r_i \le b(sk-k+1) =(sk-s+1)(sk+1)(sk-k+1), \]
hence
\begin{equation}
\label{bi}
 b_{I}\le (sk+1)(sk-k+1),
\end{equation}
and the equality $b_I =(sk+1)(sk-k+1)$ holds if and only if
\[ r_1 = r_2 = \cdots = r_b = sk - k +1. \]
Now we define an incidence structure $D^{\perp}$ on a set $X^{\perp}$
of $m$ points labeled by the $m$ compatible resolutions $R_1,\ldots, R_m$ from
$S=\{ R_1,\ldots, R_m \}$. The "blocks" of $D^{\perp}$
are labeled by the parallel classes of $D$ which appear in resolutions
from $S$, that is, by the blocks of $\cal{I}$.
By definition, a point $D^{\perp}$ labeled by $R_i$ is incident with
$r$ blocks of $D^{\perp}$ which are labeled by the $r$ blocks of $\cal{I}$
being the parallel classes of $D$ appearing in the resolution $R_i$.
Since every two resolutions $R_i$, $R_j$ ($i\neq j$) of $S$ share exactly
one parallel class, every two distinct points of $D^{\perp}$ appear together
in exactly one block of $D^{\perp}$.

Let $k_j$ ($1\le j \le b_I$) denote the size of the $j$th block of $D^{\perp}$. We have
\begin{eqnarray*}
  \sum_{j=1}^{b_I} k_j & = & mr, \\
\sum_{j=1}^{b_I} {k_j}(k_j -1) & = & m(m-1),
\end{eqnarray*}
whence\[
\sum_{j=1}^{b_I} {k^2_j}= m(m-1+r), \]
\[ ( \sum_{j=1}^{b_I} k_j )^2 = m^{2}r^2 \le {b_I}\sum_{j=1}^{b_I} k^2_j =
{b_I}m(m-1+r), \]
and
\[
mr^2 \le {b_I}(m-1+r).
\]
Applying the inequality (\ref{bi}), we have
\begin{equation}
\label{e2}
 mr^2 \le (sk+1)(sk-k+1)(m-1+r).
\end{equation}
After the substitution $r=sk+1$, inequality (\ref{e2})
simplifies to
\begin{equation}
\label{e3}
m \le (sk-k+1)s.
\end{equation}
Now assume that equality holds in (\ref{e3}), that is,
\[ m=(sk-k+1)s, \]
which is possible only if $b_I$ meets the equality in (\ref{bi}),
that is,
\[ b_I = (sk+1)(sk-k+1)= r(sk-k+1) =rn. \]
Then
\[ \sum_{j=1}^{b_I} (k_j -s)^2 = \sum_{j=1}^{b_I} {k^2_j} -
2s \sum_{j=1}^{b_I} k_j +s^2{b_I}=0, \]
thus, all blocks of $D^{\perp}$ are of size $s$, and
 $D^{\perp}$ is a 2-$((sk-k+1)s, s,1)$ design.

Suppose that $m=(sk-k+1)s$.
We define an incidence structure $\cal{P}$ with points labeled by the
$v=(sk-s+1)k$ points of $D$ and the
\[ b^\perp = b_I =(sk+1)(sk-k+1) \]
blocks of $D^{\perp}$.
Thus, $\cal{P}$ has
\[ \bar{v}=v+b^\perp =(sk-s+1)k +(sk+1)(sk-k+1) = (sk)^2 +sk +1 \]
points.

The blocks of $\cal{P}$ are of two types. The blocks of the first type are
\[ b=(s(k-1)+1)(sk+1) \]
blocks, each being a union of a block $B$ of $D$ with the block of the
dual structure of $\cal{I}$, associated with
the point of $\cal{I}$ labeled by $B$. Each such block is of size
\[ k + (sk - k+1)=sk +1 = r. \]
The blocks of the second type are also of size $r=sk+1$ and are labeled
by the $m=(sk-k+1)s$ points of $D^\perp$, and coincide with the blocks
of the dual structure of $D^\perp$.
Thus, $\cal{P}$ has
\[ \bar{b}=b+m = (sk)^2 +sk +1 \]
blocks.

Since every two points of $D^\perp$ appear in exactly one block of $D^\perp$,
every two blocks of $\cal{P}$ of the second type  intersect in exactly one 
point.
Since the $r$ blocks of $D^\perp$ through a point of $D^\perp$ correspond to
a parallel class of blocks of $\cal{I}$, every block of $\cal{P}$ of the first
type meets every block of the second type in exactly one point.

To determine how pairs of blocks of the first type intersect, we consider
the block graph $\Gamma$ of $D$, that is, the vertices of  $\Gamma$
are labeled by the blocks of $D$, where two vertices are adjacent if and only if
the corresponding blocks are not disjoint. The graph  $\Gamma$
is strongly regular with parameters
\[ (b,a,\lambda,\mu), \]
where $b=(s(k-1)+1)(sk+1)$ is the number of vertices,
\[ a=k(r-1)=sk^2 \]
is the degree of each vertex, and
\[ \lambda = (r-2) + (k-1)^2 = k(k+s -2), \ \mu =k^2. \]
The complementary graph $\bar{\Gamma}$ is strongly regular
of degree
\[ \bar{a}=b-1-a =s^{2}k^2 +2sk -s^{2}k -k^{2}s -s. \]
The $b_I =(sk+1)(sk-k+1)$ blocks of $\cal{I}$, being parallel classes
 of $D$ belonging to resolutions of $S$, are cocliques
in $\Gamma$ of size
\[ \frac{v}{k}=n=sk-s+1. \]
By the property of $S$, every two blocks of $\cal{I}$
can share at most one point, hence the corresponding $(sk-k+1)$-cocliques 
of $\Gamma$,
can share at most one vertex.
Viewed as $(sk-k+1)$-cliques of $\bar{\Gamma}$, the blocks of $\cal{I}$ cover
\[ b_{I}{ sk-k+1 \choose 2}=\frac{(sk+1)(sk-k+1)(sk-s+1)(sk-s)}{2} \]
edges of  $\bar{\Gamma}$, which is equal to the
 total number of edges of  $\bar{\Gamma}$.
It follows that every two blocks of $\cal{P}$ of the 
first type share exactly one point.
Thus, $\cal{P}$ is a projective plane of order $sk$ in which $D$ is embedded
as a maximal  $\{ (sk-s+1)k; k \}$-arc.

This completes the proof.

\begin{Rem}
\label{rem2}
{\rm
It is an interesting open question whether the results from this section can be generalized
to designs having $m$ mutually compatible resolutions, with $m$ slightly smaller than
$(sk - k+1)s$, in the spirit of the results by
A. Beutelspacher and K. Metsch on partial projective planes \cite{BM}, \cite{BM2}.
}
\end{Rem}

\section{On a conjecture by Andries Brouwer}

It is conceivable that a resolvable 2-$((sk -s +1)k, k, 1)$ design admitting
a set of mutually compatible resolutions that achieves the bound of Theorem
\ref{th}  possesses a high degree of symmetry.
One measure of symmetry is the $p$-rank of the incidence matrix
over a finite field of characteristic $p$ that divides $r-1=sk$, which is
the order of the related projective plane.

The special case $s=2, \ k=2^{t-1}, \ t\ge 2$ corresponds to projective planes 
of
order $2^t$. A 2-$(2^{2t-1}-2^{t-1}, 2^{t-1},1)$ design arising from a maximal
$(2^{2t-1}-2^{t-1}; 2^{t-1})$-arc $\cal{A}$ is called an
{\it oval} design \cite[8.4]{AK}, in reference to the fact that
the corresponding maximal $(2^t +2;2)$-arc ${\cal{A}}^{\perp}$ in the dual plane 
is a hyperoval.
The 2-rank of oval designs in $PG(2^t,2)$ has been studied extensively.
In 1989 Mackenzie \cite{M-dis}
(see also \cite[Theorem 8.4.1]{AK}, \cite{KM}) proved that the 2-rank of 
any oval design in
$PG(2^t,2)$ is bounded from above by $3^t -2^t$. It was conjectured by Assmus
that this upper bound is always achievable. This conjecture was proved
consequently by Carpenter \cite{C96}, 
by using a result by Blokhuis and Moorhouse
\cite{BM}.

In the smallest case, $t=2$, a maximal $(6;2)$-arc ($s=k=2$) in the plane of
order 4 and its dual plane is a hyperoval, yielding the unique
 trivial 2-$(6,2,1)$ design. The next case, $t=3$, corresponds 
to the  projective
plane of order 8, which contains only one class (up to projective equivalence)
of hyperovals, or $(10;2)$-arcs, and consequently, one (up to isomorphism)
maximal $(28;4)$-arc, yielding a 2-$(28,4,1)$ oval design.

Designs with the latter parameters, 2-$(28,4,1)$, have been the subject of
several papers (\cite{ABB}, \cite{BBT}, \cite{Br}, \cite{JT}, \cite{Kr}, \cite{MTW}).
According to the Handbook of Combinatorial Designs \cite[page 37]{CRC},
there are at least 4,747 known nonisomorphic designs with these parameters,
and all  designs possessing nontrivial automorphisms have been
classified (Kr\'{c}adinac \cite{Kr}). The more recent publication by
Al-Azemi,  Anton Betten, and Dieter Betten \cite{ABB} gives a much bigger number
of nonisomorphic 2-$(28,4,1)$ designs, namely, 
68,806 such designs, all having a blocking set,  and among those, 68,484 designs
have a trivial automorphism group.

In \cite{Br}, Andries Brouwer investigated the embeddability of
2-$(28,4,1)$ designs as unitals in projective planes
of order 9.
The 2-ranks of the 138 designs 
examined by Brouwer in \cite{Br} ranged between 19 and 27, with the exception of 2-rank 20. 
The minimum 2-rank,  19,
was achieved by a design being the smallest member of the family of Ree unitals, 
and one of the  two
2-$(28,4,1)$ designs having 2-transitive automorphism groups (the second being
the classical Hermitian unital, having 2-rank 21).
It was shown in \cite{JT} that there are no 2-$(28,4,1)$ designs of
2-rank 20, and there are
exactly four nonisomorphic designs of 2-rank 21, one being the classical
Hermitian unital.
  
 It turns out that the 2-$(28,4,1)$ Ree unital is isomorphic to the
oval design in the plane of order 8, $PG(2^3,2)$.
 
Brouwer \cite{Br} made the conjecture that 19 is the minimum 2-rank of any 
2-$(28,4,1)$ design, and this minimum is achieved by  the Ree unital
only. This conjecture was proved to be true by McGuire, Tonchev and Ward \cite{MTW}.

Taking into account Carpenter's result about the 2-rank
of oval designs in $PG(2^t,2)$ \cite{C96}, it is tempting
to believe that the following generalization of Brouwer's conjecture
is true.

\begin{Con}.
If $\cal{D}$ is a 2-$(2^{2t-1}-2^{t-1}, 2^{t-1},1)$ design ($t \ge 2$),
with an incidence matrix $A$, then
\[ rank_{2}(A) \ge 3^t - 2^t, \]
and the equality
\[   rank_{2}(A) = 3^t - 2^t \]
holds if and only if $\cal{D}$ is embeddable as a maximal $(2^{2t-1}-2^{t-1}; 2^{t-1})$-arc
in $PG(2^t, 2)$.
\end{Con}

The conjecture is trivially true for $t=2$, and its validity for $t=3$
follows from the results of \cite{MTW}.

\section{Acknowledgments}

The author wish to thank the anonymous referees for carefully reading the manuscript
and suggesting several improvements, including the open question formulated in Remark \ref{rem2}.

\end{document}